%% file: ex_article.tex
\documentclass[onefignum,onetabnum]{siamart171218}

\input{ex_shared}

\myexternaldocument{ex_supplement}

\usepackage{tikz}
\usetikzlibrary{datavisualization}
\usetikzlibrary{datavisualization.formats.functions}

\usepackage[group-separator={,}]{siunitx}

\newcommand{\CAN}{\mbox{$\mathsf{CAN}$}}
\newcommand{\CA}{\mbox{$\mathsf{CA}$}}

\newcommand{\xaxislabel}{Number of Columns $k$}
\newcommand{\yaxislabel}{Number of Rows $N$}
\newcommand{\xmin}{0}
\newcommand{\xmax}{2000}
\newcommand{\filename}{data/all_results.csv}
\newcommand{\twostagefile}{data/two_stage_l2_compare.csv}
\newcommand{\legendxpos}{1}
\newcommand{\legendypos}{0.35}

\begin{document}

\maketitle

\begin{abstract}
    A \emph{covering array} is an $N \times k$ array of elements from a $v$-ary alphabet such that every $N \times t$ subarray contains all $v^t$ tuples from the alphabet of size $t$ at least $\lambda$ times; this is denoted as $\CA_\lambda(N; t, k, v)$.
    Covering arrays have applications in the testing of large-scale complex systems; in systems that are nondeterministic, increasing $\lambda$ gives greater confidence in the system's correctness. 
    The \emph{covering array number}, $\CAN_\lambda(t,k,v)$ is the smallest number of rows for which a covering array on the other parameters exists.
    For general $\lambda$, only several nontrivial bounds are known, the smallest of which was asymptotically $\log k + \lambda \log \log k + o(\lambda)$ when $v, t$ are fixed.
    Additionally it has been conjectured that the $\log \log k$ term can be removed. 
    First, we affirm the conjecture by deriving an asymptotically optimal bound for $\CAN_\lambda(t,k,v)$ for general $\lambda$ and when $v, t$ are constant using the Stein--Lov\'asz--Johnson paradigm.
    Second, we improve upon the constants of this method using the Lov\'asz local lemma.
    Third, when $\lambda=2$, we extend a two-stage paradigm of Sarkar and Colbourn that improves on the general bound and often produces better bounds than even when $\lambda=1$ of other results.
    Fourth, we extend this two-stage paradigm further for general $\lambda$ to obtain an even stronger upper bound, including using graph coloring.
    And finally, we determine a bound on how large $\lambda$ can be for when the number of rows is fixed.
\end{abstract}



\section{Introduction}

Let $N, t, k, v, \lambda$ be positive integers. 
A \emph{covering array of index $\lambda$} is an $N \times k$ array $A$, where each entry is picked from a $v$-ary alphabet $\Sigma$, such that for every set of $t$ columns $S = \{s_1, \cdots, s_t\}$, the restriction of the columns of $A$ to $S$ has that each of the $v^t$ tuples $x \in \Sigma^t$ appears in these columns at least $\lambda$ times each.
We say each of these $t$-tuples is \emph{covered} if they appear at least $\lambda$ times in this way, and an \emph{interaction} to be a set of column/value tuples of size $t$; if a $t$-tuple only appears $i$ times, then we call it \emph{$i$-covered}.
Our paper uses the notation $\CA_\lambda(N; t, k, v)$ for this object.
We provide an equivalent formulation; let $[k]$ denote the set $\{1, \cdots, k\}$, and $\binom{[k]}{t}$ to be all subsets of $[k]$ of size $t$.
For a set of $t$ distinct columns $\{c_1, \cdots, c_t\} \subseteq \binom{[k]}{t}$, and a tuple of $t$ values $(v_1, \cdots, v_t) \in \Sigma^t$, define a \emph{$t$-way interaction} $I$ to be the set $\{(c_i, v_i) : 1 \le i \le t\}$. 
Let $I_{k,t}$ denote all $t$-way interactions over columns $[k]$ and values in $\Sigma$. 
Then a covering array of index $\lambda$ contains all interactions in $I_{k,t}$ at least $\lambda$ times each. 

Noteworthy applications of covering arrays include software and hardware testing \cite{colbourn_combinatorial_2004}, malware analysis \cite{leach_evolutionary_2019,ahmadi_mimosa_2021}, and machine learning \cite{cody_systematic_2022}.
Each of the rows of a covering array corresponds to a ``test'' performed. 
An important domain for covering arrays when $\lambda > 1$ is with any testing environment that is nondeterministic.
If such a scenario occurs, then running the same test multiple times may result in different outcomes.
By increasing $\lambda$, then a tester becomes more confident in the correctness of the suite of tests performed, as any fault in the system is less likely of not being detected by tests exhibiting redundancy. 

To construct a covering array with an index $\lambda > 1$, one could in principle use a covering array with index $1$ and a perfect hash family \cite{dougherty_perfect_2020} or distributing hash family \cite{colbourn_distributing_2019} of index $\lambda$.
These objects are much smaller in the number of rows compared to covering arrays, but a simple ``product'' construction produces a covering array with much higher index. 
If the hash family has index $\alpha$ and a (compatible) covering array has index $\beta$, then the resulting covering array has index $\alpha \cdot \beta$.
Covering arrays that have some interactions covered a required number of times are used as lower bounds for sequence covering arrays \cite{chee_sequence_2013}, used in event sequence testing. 
However, the number of rows (which correspond to tests) with these methods is often much larger than necessary to guarantee coverage of all interactions.

Much research has been devoted to determine the \emph{smallest} $N$ for which a covering array exists \cite{colbourn_combinatorial_2004}; we define $\CAN_\lambda(t, k, v)$ to be the minimum $N$ for which a $\CA_\lambda(N; t, k, v)$ exists.
A covering array $\CA_\lambda(N; t, k, v)$ for which $N = \CAN_\lambda(t,k,v)$ is \emph{row-optimal}.
For $\lambda = 1$, previous work has shown that $\CAN_1(t,k,v) = \Theta_{v,t}(\log k)$ \cite{colbourn_combinatorial_2004}, where the hidden constant depends on $v,t$, and is approximately $v^t$; subsequent work has attempted at improving this constant (both for the upper and lower bounds).
For $v = t = 2$, the precise value of $\CAN_1(2, k, 2)$ is known for all $k$ due to Kleitman and Spencer \cite{kleitman_families_1973}, and independently by Katona \cite{katona_two_1973}.
For specific values of $t, k, v$, several $\CAN_1(t,k,v)$ values are known \cite{kokkala_structure_2020}.
However, no other combination of $t, v, \lambda$ is $\CAN_\lambda(t,k,v)$ known for all $k$.
We investigate upper bounds on $\CAN_\lambda(t,k,v)$. 
We note that when not all $t$-sets of columns need to have this coverage property, some families of so-called ``variable-strength covering arrays'' have been shown to exhibit sub-logarithmic growth in $k$ \cite{moura_upper_2019}.

For general $\lambda$, it is the case that $\CAN_\lambda(t,k,v) = O_{v,t}(\lambda \log k)$:\footnote{Our notation $O_{v,t}(\cdots)$ indicates that $v,t$ are treated as constants.}
one can vertically juxtapose a row-optimal $\CA_1$ $\lambda$ times.
Since every interaction is covered at least once in the $\CA_1$, the vertical duplication has each interaction covered at least $\lambda$ times. 
Godbole, Skipper, and Sunley \cite{godbole_t-covering_1996} proved an asymptotically better upper bound, in that $\CAN_\lambda(t,k,v) = O(v^t t \log k + v^t \lambda \log \log k)$; they used the probabilistic method \cite{alon_probabilistic_2004}.
The bound they provided is only asymptotic, but can be easily converted into an explicit upper bound. 
More recently, Dougherty \cite{dougherty_hash_2019} showed an equivalent, also only asymptotic bound for perfect hash families, namely a bound of the form $C_1 \log k + C_2 \lambda \log \log k + o(\lambda)$, when $v,t$ are constants and $C_1, C_2$ are constants only depending on $v,t$.
The methods in that work can be applied to covering arrays; one can also derive an explicit upper bound. 
Like Godbole, Skipper, and Sunley, the proof of Dougherty also used the probabilistic method, but also employed a different technique.
Dougherty conjectured that the $o(\lambda)$ term in his bound can be removed.

In this paper, we not only affirm the conjecture, but additionally find an asymptotically optimal explicit upper bound on $\CAN_\lambda(t,k,v)$ for when $v, t$ are fixed. 
A sketch of the asymptotics we derive in this paper appears in \cite{dougherty_algorithmic_2022}, but does not provide details, analysis, an explicit bound, nor improvements to the underlying constants. 

This paper is organized as follows. 
\Cref{sec:basic_facts} proves basic facts about the covering array number $\CAN_\lambda(t,k,v)$.
\Cref{sec:lambert_w} introduces the Lambert $W$ function, which is a crucial tool for obtaining our asymptotically optimal results, along with several lemmas about $W$.
\Cref{sec:slj} provides the first asymptotically optimal upper bound on $\CAN_\lambda(t,k,v)$ (when $v, t$ are fixed) by using the $W$ function and applying otherwise standard probabilistic arguments. 
In \Cref{sec:LLL} we note that if the number of columns is sufficiently large, then the random events have some independence; we use the (symmetric) Lov\'asz local lemma to exploit that independence to obtain a tighter upper bound. 
\Cref{sec:two_stage} contains further improvements using a ``two-stage'' paradigm based on having a dedicated coverage method when the coverage gain from a randomized method becomes sufficiently small. 
In \Cref{sec:relaxation} we determine a bound on how much redundancy is possible for a given number of rows. 
And in \Cref{sec:conclusion} we conclude and provide worthwhile future research directions.
All upper bound calculations were carried out in Python using the SciPy \cite{virtanen_scipy_2020} package.

\section{Basic $\CAN$ Facts}\label{sec:basic_facts}

In this section we prove basic facts about $\CAN_\lambda(t,k,v)$, which are all easy generalizations of facts about $\CAN_1(t,k,v)$.
For more information when $\lambda=1$, see Colbourn \cite{colbourn_combinatorial_2004}.

A simple lower bound is that $\CAN_\lambda(t,k,v) \ge \lambda \cdot v^t$, as each set of $t$ columns must have at least $\lambda$ instances of each interaction, and there are $v^t$ of them.
Additionally, $\CAN_\lambda(t,k+1,v) \ge \CAN_\lambda(t,k,v)$ as every interaction within the first $k$ columns must also be $\lambda$-covered within any $\CA_\lambda(N; t, k+1, v)$.
And finally, $\CAN_\lambda(t,k,v) \le \CAN_{\lambda+1}(t,k,v)$ as any covering array with every interaction covered at least $\lambda+1$ times must automatically cover each $\lambda$ times.

Suppose $A$ is a $\CA_\lambda(N; t, k, v)$; deleting any column and any rows that contain a fixed symbol within this column yields the inequality $\CAN_\lambda(t-1, k-1, v) \le \frac{1}{v} \CAN(t,k,v)$.
The deletion of any $\lambda-1$ rows of $A$ has every interaction still covered at least once.
This implies that $\CAN_\lambda(t,k,v) \ge \CAN_1(t,k,v) + \lambda-1$.
On the other hand, vertical juxtaposition of a $\CA_1(N;t,k,v)$ $\lambda$ times has every interaction now covered at least $\lambda$ times, implying that $\CAN_\lambda(t,k,v) \le \lambda \cdot \CAN_1(t,k,v)$.

\section{The Lambert $W$ Function}\label{sec:lambert_w}

\begin{figure}
\centering
\begin{tikzpicture}
  \begin{axis}[
      xmin=-0.6,
      ymax=1,
      samples=1000,
      enlarge y limits=false,
      every axis plot/.append style={ultra thick},
      axis lines=middle
    ]
    \addplot[magenta](-0.1,x);
    \addplot [red!80!black, domain=-5:-1] (x * exp(x), x);
    \addplot [blue!80!black, domain=-1:0.1] (x * exp(x), x);
  \end{axis}

\end{tikzpicture}
  \caption{\label{fig:lambert_function}The Lambert $W$ function, with $W_0$ in blue (i.e., when $y \ge -1/e$), and $W_{-1}$ in red (i.e., when $y \le -1/e$). The magenta (vertical) line corresponds to a negative input that yields two real solutions to $W$.}
\end{figure}
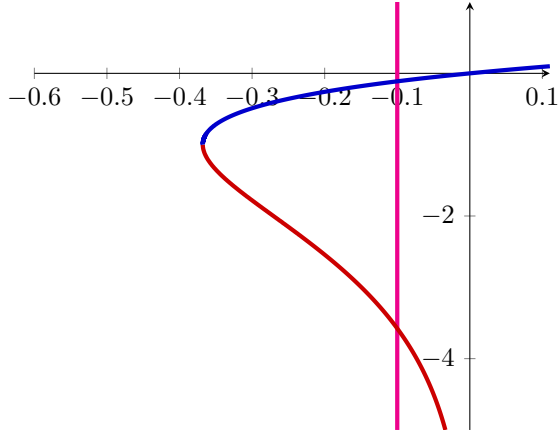

The Lambert $W$ function is defined to be the inverse of $f(W) = W \exp(W)$, and there are real solution(s) to $W(x)$ if $x \ge -1/e$.
If $-1/e < x < 0$, then $W(x)$ has two distinct real solutions, as shown in \Cref{fig:lambert_function}.
As our setup will operate in this regime of $x$, call the larger of the two solutions $W_0(x)$, and the smaller of the two $W_{-1}(x)$.

\begin{lemma}\label{lem:w_arg_is_small_enough}
    Let $t, k, v, \lambda$ be positive integers such that $k \ge t \ge 2$, $v \ge 2$, $p = 1/v^t \le 1/4$, and $a = \sqrt{\frac{(1-p)^{2\lambda} - p^{2\lambda}}{1-2p}}$.
    Denote $x = \frac{\log(1-p)}{e(\binom{k}{t} v^t a (1-p))^{1/\lambda}}$.
    Then $-1/e < x < 0$.
\end{lemma}

\begin{proof}
    $x$ is negative because $\log(1-p) < 0$, and all other quantities are positive.
    The inequality $-1/e < x$ is equivalent to:
    \[
        (1-p)^{1-1/\lambda} \left( \frac{1}{v^t a} \right)^{1/\lambda} \log \frac{1}{1-p}  < 1.
    \]
    Simple algebra and upper bounds for the left side of the above expression show this inequality is true.
\end{proof}

\begin{lemma}\label{lem:must_choose_negative_branch}
    Let $t, k, v, \lambda, p, a$ be as in \Cref{lem:w_arg_is_small_enough}.
    Then
    \[
        \frac{1}{\log(1-p)}W_0\left(\frac{\log(1-p)}{e(\binom{k}{t} v^t a (1-p))^{1/\lambda}}\right) < v^t.
    \]
\end{lemma}

\begin{proof}
    Note that $\log(1-p) < 0$, and thus the argument to $W_0(\cdot)$ is negative.
    The argument to $W_0(\cdot)$ is strictly between $-1/e$ and 0 by \Cref{lem:w_arg_is_small_enough}, and so guarantees that $W_0(x)$ is a real (negative) number.
    $W_0(-1/e) = -1$ and $-1/e$ is the only real number for which $W_0$ achieves this value (all other values $W_0$ attains are larger than $-1$).
    Since $1/\log(1-x) > -1/x$ for all $0 < x < 1$, the lemma statement can be verified with simple algebra.
\end{proof}

We now state a simple fact about $W$ that can be verified by routine calculation, and applying the definition of $W$:

\begin{lemma}\label{lem:simple_fact_about_w}
    If $ab^n n^c = d$, then $n = \frac{c}{\log b} W\left(\frac{1}{c} \left(\frac{d}{a}\right)^{1/c} \log b \right)$.
\end{lemma}


\section{Stein--Lov\'asz--Johnson Bounds}\label{sec:slj}

The methods of Stein \cite{stein_two_1974}, Lov\'asz \cite{lovasz_ratio_1975}, and Johnson \cite{johnson_approximation_1974} have been applied to covering arrays, and a proof appears in the work of Sarkar and Colbourn \cite{sarkar_upper_2017}.
We generalize their proof to provide an asymptotically tight upper bound on $\CAN_\lambda(t, k, v)$ for any positive integers $t, k, v, \lambda$.
The analysis and definitions made in the below theorem will be used throughout the rest of the paper.
We refer to this method/bound from here on as SLJ.

\begin{theorem}\label{thm:standard_prob_analysis}
	Let $t, k, v, \lambda$ be positive integers such that $k \ge t \ge 2$.
	Denote $p = 1/v^t$, and $a = \sqrt{\frac{(1-p)^{2\lambda} - p^{2\lambda}}{1-2p}}$.
	Then 
	\[
		\CAN_\lambda(t, k, v) \le 1 + \frac{\lambda}{\log(1-p)} W_{-1}\left(\frac{\log(1-p)}{e(\binom{k}{t} v^t a (1-p))^{1/\lambda}}\right).
	\]
\end{theorem}

\begin{proof}
	Let $N$ be an integer to be determined later, and let $A$ be an $N \times k$ array in which each entry is uniformly and independently selected from a $v$-ary alphabet. 
	The probability that a given interaction $T$ is not $\lambda$-covered in $A$ is $\sum_{i=0}^{\lambda-1} \binom{N}{i} p^i (1-p)^{N-i}$.
	The expected number of non-$\lambda$-covered interactions in $A$, therefore, is 
 \[
    \binom{k}{t} v^t \sum_{i=0}^{\lambda-1} \binom{N}{i} p^i (1-p)^{N-i}.
 \]
	Since for any fixed array the number of interactions not $\lambda$-covered is always an integer, if 
	\begin{equation}\label{eqn:initial_prob}
	    \binom{k}{t} v^t \sum_{i=0}^{\lambda-1} \binom{N}{i} p^i (1-p)^{N-i} < 1,
	\end{equation}
	then $A$ has positive probability of being a $\CA_\lambda$, thus proving that $\CAN_\lambda(t, k, v) \le N$.
	We repeatedly find upper bounds on the left-hand side of \Cref{eqn:initial_prob}.
	We first use the Cauchy-Schwarz inequality to give an upper bound on the summation:
	\[
		\sum_{i=0}^{\lambda-1} \binom{N}{i} p^i (1-p)^{N-i} \le \sqrt{\left(\sum_{i=0}^{\lambda-1} p^{2i} (1-p)^{2N-2i}\right) \left(\sum_{i=0}^{\lambda-1}\binom{N}{i}^2\right)}.
	\]
	The quantity $\sqrt{\sum_{i=0}^{\lambda-1} p^{2i} (1-p)^{2N-2i}}$ is equal to $(1-p)^{N-\lambda+1} \cdot a$.
	Since $\sqrt{x^2+y^2} \le x+y$ for all $x, y \ge 0$, it follows that 
	\begin{equation}
		\label{eqn:binom_inequality}
		\sqrt{\sum_{i=0}^{\lambda-1}\binom{N}{i}^2} \le \sum_{i=0}^{\lambda-1}\binom{N}{i} \le \left(\frac{eN}{\lambda}\right)^{\lambda},
	\end{equation}
	where the last inequality can be proven via induction on $\lambda$.\footnote{Note that one can improve both instances of $\lambda$ to $\lambda-1$ in the right-hand side of this inequality, which will slightly improve the constants in the derived upper bound.
However, the improvements are much smaller than what we will achieve in later sections, do not provide any significant insight, and make the theorem statement harder to understand. Therefore, we omit this from the theorem statement.}
    Consider the following equation:
	\begin{equation}\label{eqn:standard_upper_bd}
	    \binom{k}{t} v^t (1-p)^{N-\lambda+1} a \left(\frac{eN}{\lambda}\right)^{\lambda} = 1.
	\end{equation}
 	If one solves for $N$ in this equation and then adds 1 to $N$, then $N$ is an upper bound on $\CAN_\lambda$.
    Such an $N$ exists because $N^\lambda$ is a polynomial and $(1-p)^N$ is a decaying exponential; thus this equation will become less than 1 for sufficiently large $N$.
	We apply \Cref{lem:simple_fact_about_w} to obtain the following upper bound on $N$:
	\begin{equation}\label{eqn:general_result_N}
		N \le 1 + \frac{\lambda}{\log(1-p)} W\left(\frac{\log(1-p)}{e(\binom{k}{t} v^t a (1-p))^{1/\lambda}}\right).
	\end{equation}
	The argument to $W(\cdot)$ is negative since the numerator is negative and the denominator is positive; additionally, it is larger than $-1/e$, by \Cref{lem:w_arg_is_small_enough}.
	Therefore, there are two solutions $y_0, y_{-1}$ to $y_i = W(\cdot)$, where $y_{-1} \le y_0 < 0$, as is shown in \Cref{fig:lambert_function}.
	By \Cref{lem:must_choose_negative_branch}, we must choose $y_{-1}$, since if $y_0$ is chosen, then $N < \lambda v^t$, a contradiction.
	
\end{proof}

We remark that in the above proof, $a$ is very closely upper bounded by $(1-p)^\lambda \sqrt{\frac{1}{1-2p}}$, which cancels out the $(1-p)^{-\lambda}$ term in \Cref{eqn:standard_upper_bd}. 
Therefore, the only ``strong'' dependence on $\lambda$ is with the $\left(\frac{eN}{\lambda}\right)^\lambda$ term. 

\begin{figure}
    \centering
    \includegraphics[width=\textwidth]{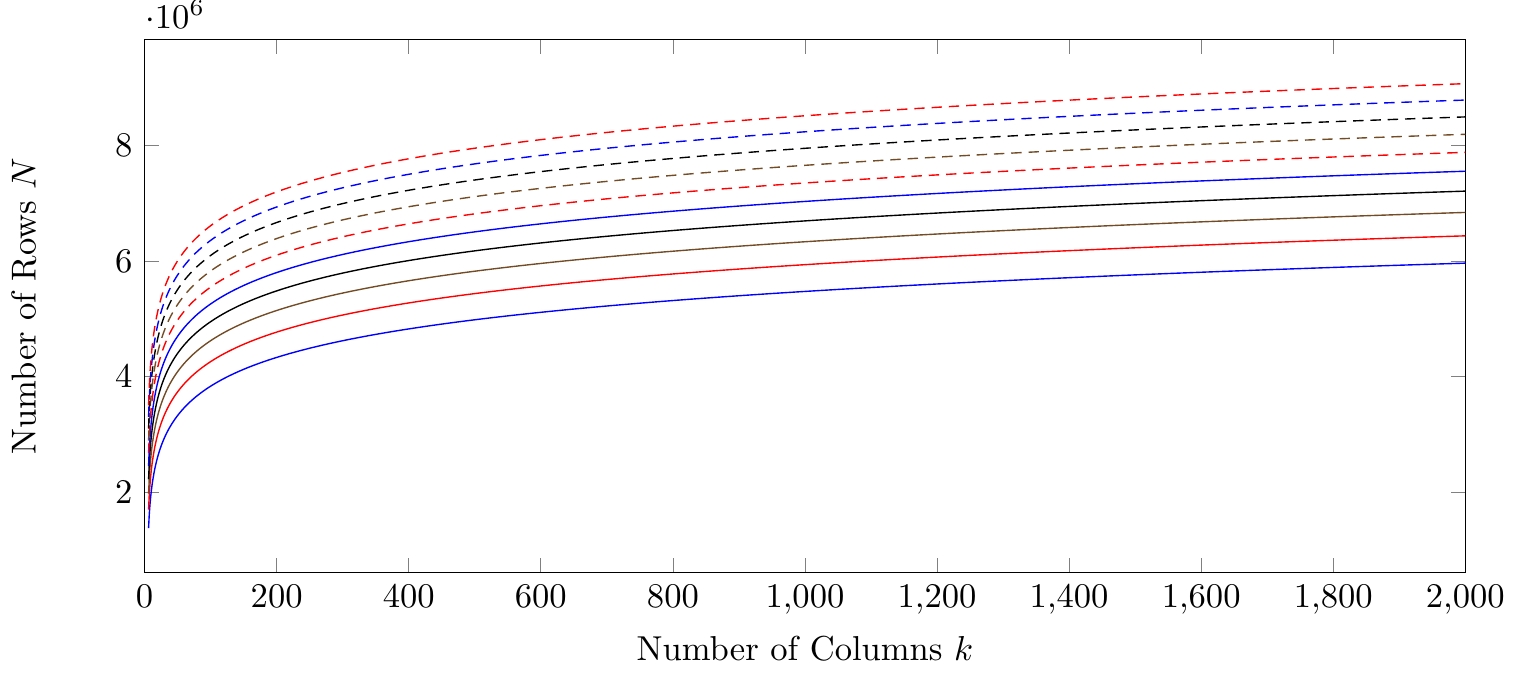}
    \caption{The upper bounds for $\CAN_\lambda(6,k,7)$ from \Cref{thm:standard_prob_analysis} for $6 \le k \le 2000$ and $1 \le \lambda \le 10$.
    The lowest curve is $\lambda=1$, and the highest is $\lambda=10$.
    }
    \label{fig:compare_lambda_values_slj}
\end{figure}

\Cref{fig:compare_lambda_values_slj} gives an example of how few rows are needed to achieve higher redundancies $\lambda$.
We plot the upper bound from \Cref{thm:standard_prob_analysis} for $\CAN_\lambda(t,k,v)$ where $t = 6, k \le 2000, v=7$, and $\lambda \le 10$.
For example, $\CAN_1(6,2000,7) \le \num{5964087}$ and $\CAN_{10}(6,2000,7) \le \num{9073425}$. 
With much fewer than twice the original number of rows, the amount of guaranteed redundancy is 10 times as large. 
The implications of how the argument of $W_{-1}$ impacts $\CAN_\lambda$ may not be immediately evident. 
We now provide a corollary that removes the dependence on $W_{-1}$, with a slightly worse bound than \Cref{thm:standard_prob_analysis} has.

\begin{corollary}\label{cor:phrase_first_bound_without_lambert}
        Let $t, k, v, \lambda, p, a$ be as in \Cref{thm:standard_prob_analysis}.
	Then 
	\[
		\CAN_\lambda(t, k, v) \le 1 + \frac{\lambda e}{(e-1)\log\frac{1}{1-p}} \left(1 +  \log\left(1 + \frac{(\binom{k}{t} v^t a (1-p))^{1/\lambda}}{\log\frac{1}{1-p}} \right) \right).
	\]
\end{corollary}

\begin{proof}
    Alzahrani and Salem \cite{alzahrani_sharp_2018} show that $W_{-1}(-e^{-z-1}) > -\alpha(z+1)$, where $\alpha = e/(e-1)$, and $z \ge 0$.
	Solving for $z$ using the argument from \Cref{eqn:general_result_N} in the proof of \Cref{thm:standard_prob_analysis} yields:
	\[
		z = \log\left(\frac{(\binom{k}{t} v^t a (1-p))^{1/\lambda}}{\log \frac{1}{1-p}} \right).
	\] 
	Substitution of the inequality and $z$ into \Cref{eqn:general_result_N} yields the corollary statement.
\end{proof}

In \Cref{fig:slj_bounds_v4_t6_l12} we plot upper bounds reported by \Cref{thm:standard_prob_analysis} and \Cref{cor:phrase_first_bound_without_lambert} for $\CAN_\lambda(t,k,v)$ where $t = 6, k \le 2000, v = 4$, and $\lambda = 12$.
Additionally for each $k$ we plot the minimum $N$ for which \Cref{eqn:initial_prob} is satisfied. 
As expected, \Cref{eqn:initial_prob} provides the smallest bound among the three, showcasing the price undertaken by taking the upper bounds to the sum.
However, the multiplicative difference between the three bounds will always be at most a constant independent of $k$ and $\lambda$, as we implicitly prove next. 

\begin{figure}
    \centering
    \begin{tikzpicture}
\begin{axis}[no marks,
	legend pos=south east,
        legend style={at={(\legendxpos,\legendypos)},anchor=south east},
	width=0.9\textwidth,
        height=3in,
	xmin=\xmin,
	xmax=\xmax,
	xlabel={\xaxislabel},
	ylabel={\yaxislabel}]
\addplot[black] table [mark=none, x=k, y=slj_no_sum_no_W, col sep=comma] {\filename};
\addplot[blue] table [mark=none, x=k, y=slj_no_sum_with_W, col sep=comma] {\filename};
\addplot[red] table [mark=none, x=k, y=slj_with_sum, col sep=comma] {\filename};

\legend{\Cref{cor:phrase_first_bound_without_lambert}, \Cref{thm:standard_prob_analysis}, \Cref{eqn:initial_prob}}
\end{axis}
\end{tikzpicture}
    \caption{Upper bounds for $\CAN_\lambda(t,k,v)$ when $t = 6, k \le 2000, v = 4$, and $\lambda = 12$ from \Cref{cor:phrase_first_bound_without_lambert}, \Cref{thm:standard_prob_analysis}, and \Cref{eqn:initial_prob}.}
    \label{fig:slj_bounds_v4_t6_l12}
\end{figure}

We now prove that the bound obtained in \Cref{cor:phrase_first_bound_without_lambert} is asymptotically optimal for when $v,t$ are constants, thus showing that \Cref{thm:standard_prob_analysis} is also asymptotically optimal also when $v,t$ are constants.

\begin{corollary}\label{cor:ca_tight_bound_first}
    $\CAN_\lambda(t,k,v) = \Theta_{v,t}(\log k + \lambda)$.
\end{corollary}

\begin{proof}
    The upper bound is a result of \Cref{cor:phrase_first_bound_without_lambert}.
    For the lower bound, we can assume without loss of generality that any covering array does not have two identical columns as $v, t \ge 2$.
    Then any row-optimal $\CA_1$ must have at least $\Omega_{v,t}(\log k)$ rows.
    To complete such an array to be a $\CA_\lambda$, we require at least $\lambda-1$ more rows, as some interaction is 1-covered in this array.
    Therefore, $\CAN_\lambda(t,k,v) = \Omega_{v,t}(\log k + \lambda)$.
\end{proof}

For the upper bound, an analysis of \Cref{thm:standard_prob_analysis,cor:phrase_first_bound_without_lambert} shows that the upper bound is approximately equal to 
\[
    \lambda v^t  + v^t \log \binom{k}{t}  + v^t \log v^t + v^t \log a + v^t \log (1-p) + \lambda v^t \log v^t.
\]
This bound uses the fact that $\frac{1}{\log \frac{1}{1-p}} \le v^t$ \cite{sarkar_upper_2017}. 
Since $v^t$ is much larger than $\log v^t$, a short calculation yields $\CAN_\lambda(t,k,v) = O(v^t \log \binom{k}{t} + \lambda v^t \log v^t)$.
A natural question is whether the constants involving $v, t$ can be improved. 
When $k$ is sufficiently small, then $\lambda v^t$ rows are required for sufficiently large $\lambda$; such arrays are called \emph{orthogonal arrays}, wherein each interaction is covered the same number of times ($\lambda$).
Therefore, the constant for $\lambda$ can (at best) be lowered to $v^t$ for all $k$.

The bounds in \Cref{thm:standard_prob_analysis,cor:phrase_first_bound_without_lambert} are non-constructive; however, we outline the method of Dougherty et al. \cite{dougherty_algorithmic_2022} that produces a $\CA_\lambda(N; t, k, v)$ that runs in polynomial time in $k$ (when $v, t$ are fixed) and the produced array size is at most these bounds.
This method generalizes the ``density'' method of Bryce and Colbourn \cite{bryce_density-based_2009} for covering arrays with $\lambda=1$ as their method does not easily generalize to $\lambda > 1$.
Compute the number $N$ based on $t, k, v, \lambda$ from \Cref{thm:standard_prob_analysis}. 
Build the $\CA$ one row at a time; suppose $r \ge 0$ rows have been completed, the $(r+1)$-st row is not complete, and there exists some interaction $I$ that has not been $\lambda$-covered yet in these $r$ rows.
Determine the probability that, when fixing one more entry in column $c$ of this $(r+1)$-st row to a certain value $x$, $I$ will be $\lambda$-covered in the remaining $N-r$ rows if their entries are chosen uniformly at random, including the remaining entries in the row we are building.
Now sum this probability across all interactions that are not yet $\lambda$-covered.
This is the expected number of remaining not-$\lambda$-covered interactions after all rows have been built. 
Pick any value $x$ in column $c$ such that this expectation is minimized.
Since the expectation from \Cref{thm:standard_prob_analysis} starts at a value strictly less than 1, and never increases at any point, a $\CA$ will be built that meets this bound.
A disadvantage of the density algorithm is that a counter needs to be stored for how many times each interaction has been covered so far, and there are $\binom{k}{t} v^t$ interactions.

\section{Lov\'asz local lemma Bounds}\label{sec:LLL}

The (symmetric) Lov\'asz local lemma is a well known result that involves exploiting independence among events in a probability space if such independence exists to guarantee the existence of all events simultaneously not occurring.
The formal statement of the lemma is provided next.

\begin{lemma}[\cite{lovasz_ratio_1975}]
\label{lem:LLL_sym}
    Let $A_0, \ldots, A_{m-1}$ be events in an arbitrary probability space. Suppose that for each event $A_i$, $Pr(A_i) \le p$, and that $A_i$ is independent of all other events except at most $d$ of them. If 
    \[
        ep(d+1)\le 1,
    \]
    then with non-zero probability none of the events occur.
\end{lemma}

In our context, an event is a set of $t$ columns not having all $v^t$ interactions covered.
Therefore, if $e(d+1)p\le 1$, with the probability $p$ as defined before with $N$ rows, then there exists a covering array with $N$ rows.
We now apply the local lemma to find a bound on $\CAN_\lambda$, with setup similar to that of Godbole, Skipper, and Sunley \cite{godbole_t-covering_1996}.
We refer to this method/bound from here on as LLL.

\begin{theorem}\label{thm:standard_lll}
    Let $t, k, v, \lambda, p, a$ be as in \Cref{thm:standard_prob_analysis}.
    Then 
	\[
		\CAN_\lambda(t, k, v) \le \frac{\lambda}{\log(1-p)} W_{-1}\left(\frac{\log(1-p)}{e((\binom{k}{t} - \binom{k-t}{t}) v^t a (1-p) e)^{1/\lambda}}\right).
	\]
\end{theorem}

\begin{proof}
    We proceed similarly as in \Cref{thm:standard_prob_analysis}.
    Let $A$ be an array with $N$ rows and $k$ columns that has entries from a $v$-ary alphabet chosen uniformly at random and independent of other entries.
    The probability that all $v^t$ interactions are not covered in a given set of $t$ columns of $A$ is $v^t \sum_{i=0}^{\lambda-1} \binom{N}{i} (\frac{1}{v^t})^i (1-\frac{1}{v^t})^{N-i}$.
    
    As defined previously, we associate for each set of $t$ columns an event $A_i$; $A_i$ occurs precisely when all $v^t$ interactions in these columns are covered.
    Under this definition, $A_i$ is independent of a different event $A_j$ precisely when these two events do not share a column.
    Therefore, the number of other events $A_i$ is dependent upon is precisely $d = \binom{k}{t} - \binom{k-t}{t} - 1$. 
    By \Cref{lem:LLL_sym}, if 
    \begin{equation}\label{eqn:lll_standard}
        e\left(\binom{k}{t} - \binom{k-t}{t}\right) v^t \sum_{i=0}^{\lambda-1} \binom{N}{i} \left(\frac{1}{v^t}\right)^i \left(1-\frac{1}{v^t}\right)^{N-i} \le 1,
    \end{equation}
    then a $\CA_\lambda(N; t, k, v)$ exists.
    We then apply the same method as in \Cref{thm:standard_prob_analysis} to find an upper bound on $N$ for which this equation is satisfied.
    We do not have to add 1 to the result as we only require that the above expression is less than or equal to 1, not strictly less than 1.
\end{proof}

\Cref{thm:standard_lll} improves on \Cref{thm:standard_prob_analysis} when $e\left(\binom{k}{t} - \binom{k-t}{t}\right) < \binom{k}{t}$. 
Stirling's formula for a factorial gives an asymptotic for the binomial coefficient $\binom{n}{k}$ when $n, k$ are sufficiently large.
A short calculation shows that if $k$ is sufficiently larger than $t$, then the above inequality is satisfied.
As a result of \Cref{thm:standard_lll} we achieve a very easy corollary, again using the same lower bound on $W_{-1}$ as was used in the proof of \Cref{cor:phrase_first_bound_without_lambert}.
\begin{corollary}\label{cor:lll_without_lambert}
    Let $t, k, v, \lambda, p, a$ be as in \Cref{thm:standard_prob_analysis}.
    Then 
    \[
        \CAN_\lambda(t, k, v) \le \frac{\lambda e}{(e-1)\log \frac{1}{1-p}} \left(1 +  \log\left(1 + \frac{((\binom{k}{t} - \binom{k-t}{t}) v^t a (1-p) e)^{1/\lambda}}{\log \frac{1}{1-p}} \right) \right).
    \]
\end{corollary}

\begin{figure}
    \centering
    \begin{tikzpicture}
\begin{axis}[no marks,
	legend pos=south east,
        legend style={at={(\legendxpos,\legendypos)},anchor=south east},
	width=0.9\textwidth,
        height=3in,
	xmin=\xmin,
	xmax=\xmax,
	xlabel={\xaxislabel},
	ylabel={\yaxislabel}]
\addplot[black,forget plot] table [mark=none, x=k, y=slj_no_sum_no_W, col sep=comma] {\filename};
\addplot[black,forget plot] table [mark=none, x=k, y=slj_no_sum_with_W, col sep=comma] {\filename};
\addplot[black,forget plot] table [mark=none, x=k, y=slj_with_sum, col sep=comma] {\filename};

\addplot[red] table [mark=none, x=k, y=lll_no_sum_no_W, col sep=comma] {\filename};
\addplot[blue] table [mark=none, x=k, y=lll_no_sum_with_W, col sep=comma] {\filename};
\addplot[orange] table [mark=none, x=k, y=lll_with_sum, col sep=comma] {\filename};

\legend{\Cref{cor:lll_without_lambert}, \Cref{thm:standard_lll}, \Cref{eqn:lll_standard}}
\end{axis}
\end{tikzpicture}
    \caption{Upper bounds for $\CAN_\lambda(t,k,v)$ when $t = 6, k \le 2000, v = 4$, and $\lambda = 12$ from \Cref{cor:lll_without_lambert}, \Cref{thm:standard_lll}, and \Cref{eqn:lll_standard}.
    These three bounds are the lowest lines within each of the three sections.
    The other upper bounds are from the three SLJ results from \Cref{fig:slj_bounds_v4_t6_l12}.}
    \label{fig:lll_bounds_v4_t6_l12}
\end{figure}
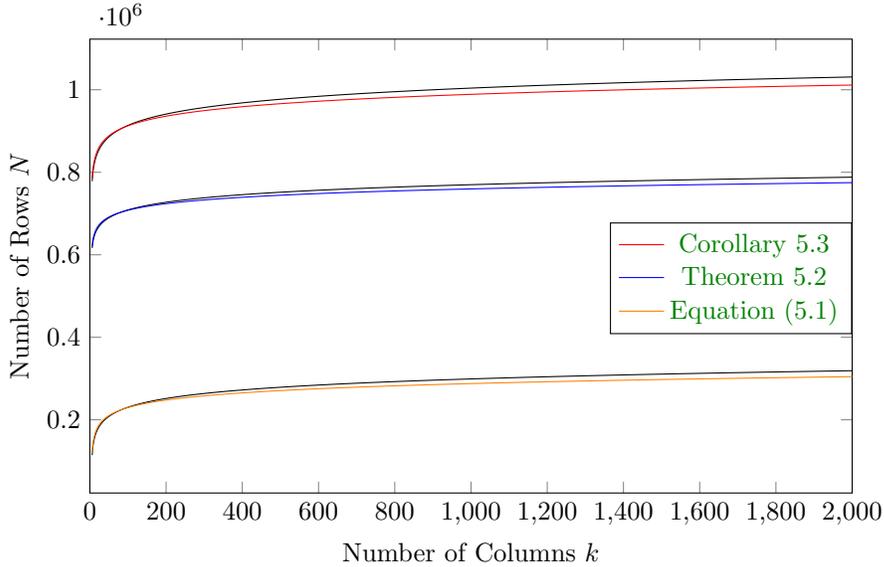

We compare the results from \Cref{cor:lll_without_lambert}, \Cref{thm:standard_lll}, and \Cref{eqn:lll_standard} in \Cref{fig:lll_bounds_v4_t6_l12} by plotting the upper bounds of both for $\CAN_\lambda(t,k,v)$ where $t = 6$, $k \le 2000$, $v = 7$, and $\lambda = 12$.
Additionally, we plot the three SLJ results from \Cref{fig:slj_bounds_v4_t6_l12} in black for comparison. 
When $k$ is sufficiently small, the upper bound from \Cref{cor:phrase_first_bound_without_lambert} is smaller than that of \Cref{cor:lll_without_lambert}; but when $k \ge 85$, then \Cref{cor:lll_without_lambert} produces a smaller upper bound.
Similar comparisons can be made from the other two bounds between LLL and SLJ.

The bounds stated in \Cref{thm:standard_lll,cor:lll_without_lambert} are non-constructive, but the methods of Moser and Tard\"os \cite{moser_constructive_2010} can make them constructive through the means of an algorithm that runs in polynomial expected time in the number of events.
We outline their method within the setting of covering arrays. 
Determine the smallest value of $N$ for which \Cref{eqn:lll_standard} is satisfied, and construct an array with $N$ rows (and $k$ columns) with entries chosen uniformly at random, independent of other entries. 
For each set of $t$ columns (equivalently, an event), check if all interactions in these columns are $\lambda$-covered; if not, resample these $t$ columns (across all rows) by choosing entries again uniformly at random, independently. 
Then repeat the check from the first step. 
The advantage of the local lemma is that, in addition to the number of rows being smaller, the coverage of each interaction only needs to be stored throughout the check of each set of $t$ columns.
Compare this to the storage requirements for the density method outlined previously; there, all interactions need to have their coverage stored at the start of the algorithm.

\section{Alteration/Two-Stage Bounds}\label{sec:two_stage}

We can improve upon the upper bounds of \Cref{thm:standard_prob_analysis,thm:standard_lll} as follows.
Create a random $N \times k$ array as before, and suppose that the number of interactions not $\lambda$-covered is at most some quantity $\rho$; we denote this as the \emph{first stage}.
Then one can add $\lambda\cdot\rho$ rows, each of which covers each of these $\rho$ interactions; this is the \emph{second stage}.
The key is determining how many initial rows $N$ to choose for the initial random array.
If $N$ is too large, then this method becomes mostly irrelevant as one can use \Cref{thm:standard_prob_analysis}.
If $N$ is too small, then this method is much less competitive because the number of rows will be much larger than a constant times $\log k + \lambda$.
Sarkar and Colbourn \cite{sarkar_upper_2017} found that when $\lambda = 1$, the optimal choice of $N$ is
\[
\frac{\log\binom{k}{t} + t \log v + \log \log \frac{1}{1-p}}{\log \frac{1}{1-p} }
\]
where $p = \frac{1}{v^t}$.
This technique is known as ``alteration'' \cite{alon_probabilistic_2004}, and is useful for constructing covering arrays; there are two reasons for why two stage approaches are useful.
The first is that ``stronger'' methods can be used in the second stage as the number of interactions left uncovered at this point is often very small, which decrease computation time. 
The number of interactions becoming $\lambda$-covered decays towards the end of the density algorithm (and other methods as well); therefore, more dedicated methods are desired when this decay becomes sufficiently small. 
The second is that such second-stage methods know what algorithm came in the first stage, and thus can in principle have more information about what interactions are left uncovered.

In the general case of arbitrary $\lambda > 2$, we cannot determine the optimal choice of the first stage's number of rows. 
However, when $\lambda = 2$, we can do so, and thus gain an improvement for $\CAN_2(t,k,v)$, as we prove next.

\begin{theorem}\label{thm:differentiation_2}
    Let $t, k, v, p$ be as in \Cref{thm:standard_prob_analysis}.
	Then
	\[
		\CAN_2(t, k, v) \le m + 2\binom{k}{t} v^t \left((1-p)^m + m p (1-p)^{m-1}\right),
	\]
        where
        \begin{equation}\label{eqn:diff_lambda2_minN}
            m = \frac{1}{\log(1-p)} \left(W_{-1}\left(\frac{-e(1-p)^{v^t}}{2\binom{k}{t}}\right) - (v^t-1)\log(1-p)-1 \right).
        \end{equation}
\end{theorem}

\begin{proof}
	The proof structure is very similar to \Cref{thm:standard_prob_analysis}.
	Let $A$ be an array with $N$ rows, with entries chosen uniformly at random, independently.
	The expected number of uncovered interactions in $A$ is $\binom{k}{t} v^t \left((1-p)^N + N p (1-p)^{N-1}\right)$.
	To complete the remaining interactions, one can add $\lambda = 2$ rows for each uncovered interaction, yielding a covering array of index $\lambda$ with
 \begin{equation}\label{eqn:diff_2_intermediate_eqn}
     N + 2 \binom{k}{t} v^t \left((1-p)^N + N p (1-p)^{N-1}\right)
 \end{equation} 
 rows. The minimum occurs when the derivative of \Cref{eqn:diff_2_intermediate_eqn} with respect to $N$ is equal to 0.
 Solving this equation for $N$ yields the minimum choice for $N$ occurs at the $m$ value specified above in \Cref{eqn:diff_lambda2_minN}.
 Substitute this value of $N$ into \Cref{eqn:diff_2_intermediate_eqn}. Note that the argument to $W_{-1}$ is always strictly between $-1/e$ and 0, thus the arguments in the proof of \Cref{thm:standard_prob_analysis} apply.
 
	
\end{proof}

\begin{figure}
    \centering
    \begin{tikzpicture}
\begin{axis}[no marks,
	legend pos=south east,
	width=0.9\textwidth,
        height=3in,
	xmin=\xmin,
	xmax=\xmax,
	xlabel={\xaxislabel},
	ylabel={\yaxislabel}]
\addplot[black] table [mark=none, x=k, y=slj_no_sum_W, col sep=comma] {\twostagefile};
\addplot[blue] table [mark=none, x=k, y=slj_sum, col sep=comma] {\twostagefile};
\addplot[red] table [mark=none, x=k, y=two_slj_sum, col sep=comma] {\twostagefile};

\legend{\Cref{thm:standard_prob_analysis}, \Cref{eqn:standard_upper_bd},  \Cref{thm:differentiation_2}}
\end{axis}
\end{tikzpicture}
    \caption{Upper bounds for $\CAN_\lambda(t,k,v)$ when $t = 6, k \le 2000, v = 4$, and $\lambda = 2$ from \Cref{eqn:standard_upper_bd}, \Cref{thm:standard_prob_analysis}, and \Cref{thm:differentiation_2}.}
    \label{fig:two_stage_comparison_l2}
\end{figure}
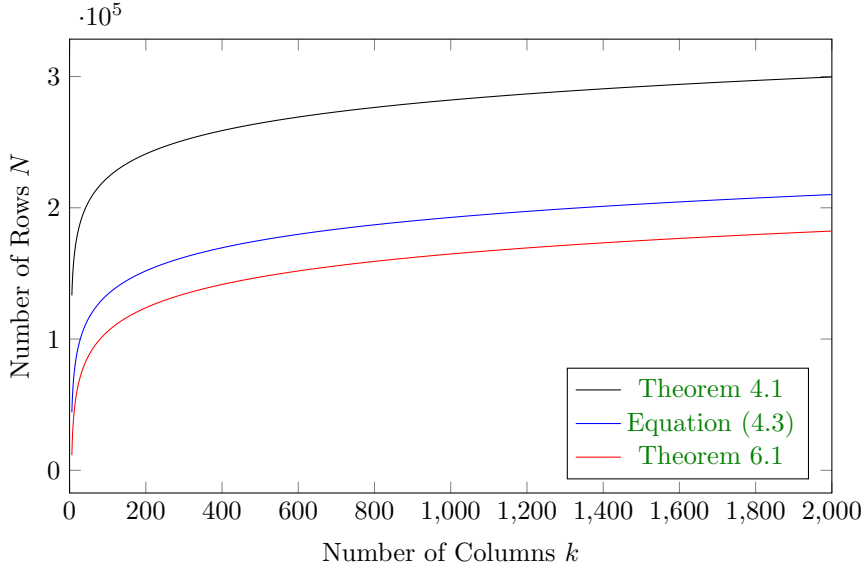

\begin{corollary}\label{cor:lambda_2_result}
    Let $t, k, v, p$ be as in \Cref{thm:standard_prob_analysis}.
	Then
	\[
		\CAN_2(t, k, v) \le m + 2\binom{k}{t} v^t \left((1-p)^m + m p (1-p)^{m-1}\right),
	\]
        where
        \[
            m = \frac{e}{e-1} \frac{\log\binom{k}{t} + v^t \log\frac{1}{1-p} + \log 2}{\log\left(\frac{1}{1-p} \right)} + 1 - v^t.
        \]
\end{corollary}

\begin{proof}
    Use the lower bound of $W_{-1}$ by Alzharani and Salem from the proof of \Cref{thm:standard_prob_analysis} within \Cref{eqn:diff_lambda2_minN}.
\end{proof}

Routine verification shows that this alteration approach improves the SLJ upper bounds when $\lambda = 2$; see \Cref{fig:two_stage_comparison_l2} for a comparison for when $t = 6, k \le 2000$, and $v=4$.
Additionally, the SLJ and LLL bounds for standard covering arrays obtained by Sarkar and Colbourn \cite{sarkar_two-stage_2019} for $\lambda = 1$ are very similar in size to that of \Cref{thm:differentiation_2,cor:lambda_2_result}.
Remarkably, \Cref{thm:differentiation_2} produces a smaller value for $\CAN_2(t,k,v)$ than their methods do for $\lambda=1$, sufficiently small $k$, and all $t,v$; \Cref{cor:lambda_2_result} gives the same results but only for smaller $k$.
The reason for this improvement is simple, although we only present why small values of $k$ are improved from \Cref{cor:lambda_2_result} for ease of analysis.
The SLJ result for $\lambda=1$ by Sarkar and Colbourn \cite{sarkar_upper_2017} proves that $\CAN_1(t,k,v)$ is approximately $v^t \log \binom{k}{t} + v^t t \log v$.
Compare that with the result of \Cref{cor:lambda_2_result} for $\lambda = 2$; after a short calculation, the produced number of rows with the chosen value of $m$ is approximately 
\[
    \frac{e}{e-1} v^t \log \binom{k}{t} \left(1 + 2 \binom{k}{t}^{-1/(e-1)} \right).
\]
When $k$ is sufficiently small, this is smaller than the SLJ result, mainly because of the constant on the rightmost term being larger for SLJ.
If $k$ is sufficiently large, the leading constant from the SLJ result yields a smaller upper bound than for \Cref{cor:lambda_2_result}.

Even though the SLJ and LLL results for $\lambda=1$ are asymptotically smaller than that of \Cref{cor:lambda_2_result} when $\lambda=2$ for sufficiently large $k$, how large does $k$ have to be for the asymptotics to take over?
The answer is often larger than is needed for practical use. 
We give an example with $t = 6, k = 2000, v = 7$.
SLJ reports that $\CAN_1(t, k, v) \le \num{5608361}$, whereas \Cref{cor:lambda_2_result} reports that $\CAN_2(t,k,v) \le \num{5236206}$; not only is this an improvement of over \num{370000} rows, but additional redundancy is guaranteed.
We demonstrate this further in \Cref{tbl:slj_vs_lll_vs_2stage_lambda2} by calculating their SLJ and LLL bounds, and the results from \Cref{thm:differentiation_2,cor:lambda_2_result} for $\lambda=2$, for all $k$ powers of 10 between $10^1$ and $10^{10}$.
As recorded in the table, for all $k \le 10^4$, \Cref{thm:differentiation_2} improves upon existing results; further, for all $k \le 10^{10}$, \Cref{thm:differentiation_2} improves upon SLJ.
More specifically, an improvement occurs for this choice of $v, t$ for all $k \le \num{34215}$.

\begin{table}[]
\begin{tabular}{c|ll|ll|}
\cline{2-5}
                          & \multicolumn{2}{c|}{$\lambda=1$}                           & \multicolumn{2}{c|}{$\lambda=2$} \\ \hline
\multicolumn{1}{|c|}{$k$} & \multicolumn{1}{c|}{SLJ} & \multicolumn{1}{c|}{LLL} & \multicolumn{1}{l|}{\Cref{thm:differentiation_2}}   & \Cref{cor:lambda_2_result}  \\ \hline
\multicolumn{1}{|l|}{$10^{1}$}    & \multicolumn{1}{l|}{\num{2002680}}                  & \num{2120329}                  & \multicolumn{1}{l|}{\textbf{\num{1089371}}}    &  \num{1214439}     \\ \hline
\multicolumn{1}{|l|}{$10^{2}$}    & \multicolumn{1}{l|}{\num{3832330}}                  &   \num{3814804}                & \multicolumn{1}{l|}{\textbf{\num{3040435}}}    &  \num{4087136}     \\ \hline
\multicolumn{1}{|l|}{$10^{3}$}    & \multicolumn{1}{l|}{\num{5473916}}                  &   \num{5199000}                & \multicolumn{1}{l|}{\textbf{\num{4734170}}}    &  \num{6684079}     \\ \hline
\multicolumn{1}{|l|}{$10^{4}$}    & \multicolumn{1}{l|}{\num{7100882}}                  &   \num{6556396}                & \multicolumn{1}{l|}{\textbf{\num{6396559}}}    &  \num{9257901}     \\ \hline
\multicolumn{1}{|l|}{$10^{5}$}    & \multicolumn{1}{l|}{\num{8726415}}                  &    \textbf{\num{7911166}}               & \multicolumn{1}{l|}{\num{8049136}}    &  \num{11829456}     \\ \hline
\multicolumn{1}{|l|}{$10^{6}$}    & \multicolumn{1}{l|}{\num{10351805}}                  &   \textbf{\num{9265673}}                & \multicolumn{1}{l|}{\num{9696435}}    &  \num{14400785}     \\ \hline
\multicolumn{1}{|l|}{$10^{7}$}    & \multicolumn{1}{l|}{\num{11977180}}                  &  \textbf{\num{10620155}}                 & \multicolumn{1}{l|}{\num{11340237}}    & \num{16972092}      \\ \hline
\multicolumn{1}{|l|}{$10^{8}$}    & \multicolumn{1}{l|}{\num{13602555}}                  &  \textbf{\num{11974633}}                 & \multicolumn{1}{l|}{\num{12981515}}    & \num{19543396}      \\ \hline
\multicolumn{1}{|l|}{$10^{9}$}    & \multicolumn{1}{l|}{\num{15227929}}                  &   \textbf{\num{13329112}}                & \multicolumn{1}{l|}{\num{14620881}}    &  \num{22114700}     \\ \hline
\multicolumn{1}{|l|}{$10^{10}$}    & \multicolumn{1}{l|}{\num{16853303}}                  &    \textbf{\num{14683590}}               & \multicolumn{1}{l|}{\num{16258748}}    &  \num{24686004}     \\ \hline
\end{tabular}
\centering
\caption{\label{tbl:slj_vs_lll_vs_2stage_lambda2}SLJ and LLL results from Sarkar and Colbourn \cite{sarkar_upper_2017} for $\lambda=1$ compared to \Cref{thm:differentiation_2} and \Cref{cor:lambda_2_result} with $\lambda=2$. 
Here, $t=6, v=7$, and $k$ all powers of 10 between $10^1$ and $10^{10}$ are presented.}
\end{table}

\subsection{Arbitrary $\lambda$}

One can attempt to generalize the methods of \Cref{thm:differentiation_2} to arbitrary $\lambda$ as follows.
The number of rows needed to cover all interactions would involve calculating the expected number of uncovered interactions and then adding $\lambda$ rows for each. 
This yields the following expression:
\begin{equation}\label{eqn:generalize_diff_to_arbitrary_lambda}
N + \lambda \binom{k}{t} v^t \sum_{i=0}^{\lambda-1} \binom{N}{i} p^i (1-p)^{N-i}.
\end{equation}
If $\lambda \ge 3$, determining when the derivative of \Cref{eqn:generalize_diff_to_arbitrary_lambda} with respect to $N$ is equal to 0 is not analytically solvable, at least in terms of the $W$ function. 
However, since this equation is convex in $N$ for $N > 0$, there exists a unique minimum and can be found computationally via any root-finding algorithm.  
Additionally, to find a general bound, one can take an upper bound on the derivative of \Cref{eqn:generalize_diff_to_arbitrary_lambda}, which will yield an upper bound on $\CAN_\lambda(t,k,v)$.
This is possible because the above equation is convex in $N$ when $N$ is positive.
Our next result summarizes an improvement over \Cref{thm:standard_prob_analysis,thm:standard_lll} using the same alteration approach as \Cref{thm:differentiation_2}.
The result also has an analogous version with $W_{-1}$, which gives stronger upper bounds, but we omit its statement for brevity and instead report its results in \Cref{fig:two_stage_bounds_v4_t6_l12} below.



\begin{theorem}\label{thm:differentiation_general}
    Let $t, k, v, \lambda, p, a$ be as in \Cref{thm:standard_prob_analysis}. 
    Then
    \[
   	\CAN_\lambda(t,k,v) \le m + \lambda \binom{k}{t} v^t \sum_{i=0}^{\lambda-1} \binom{m}{i} p^i (1-p)^{m-i},
    \]
    where
    \[
        m = 1 + \frac{\lambda e}{(e-1)\log \frac{1}{1-p}} \left(1 +  \log\left(1 + \frac{(\binom{k}{t} v^t a (1-p))^{1/\lambda}}{(\log \frac{1}{1-p})^{1-1/\lambda}} \right) \right).
    \]
\end{theorem}

\begin{proof}
    We proceed with the same setup as \Cref{thm:differentiation_2}.
    Let $A$ be an array with $N$ rows, with entries chosen uniformly at random, independent of other entries.
	The expected number of uncovered interactions in $A$ is $\binom{k}{t} v^t \sum_{i=0}^{\lambda-1} \binom{N}{i} p^i (1-p)^{N-i}$.
	To complete the remaining interactions, as in \Cref{thm:differentiation_2}, we add $\lambda$ rows for each uncovered interaction, yielding a $\CA_\lambda(N';t,k,v)$ with:
	\begin{equation}\label{eqn:diff_general_setup_original}
	    N' = N + \lambda \binom{k}{t} v^t \sum_{i=0}^{\lambda-1} \binom{N}{i} p^i (1-p)^{N-i}.
	\end{equation}
	For $N > 0$ and all other parameters fixed, \Cref{eqn:diff_general_setup_original} is convex, which implies there is a minimum $N^\star$ on the interval.
	However, since we want a bound on $N$, we will be taking upper bounds; any deviation on $N$ from $N^\star$ will still be an upper bound on $\CAN_\lambda(t,k,v)$.
	We differentiate the above expression with respect to $N$, and attempt to determine the value of $N$ for which the derivative is 0:
	\begin{equation}\label{eqn:diff_general_setup}
	    1 + \lambda \binom{k}{t} v^t \sum_{i=0}^{\lambda-1} \binom{N}{i} p^i (1-p)^{N-i} (H_N - H_{N-i} + \log(1-p)) = 0.
	\end{equation}
	Here, $H_k$ is the $k$th harmonic number, i.e., $H_k = \sum_{i=1}^{k} \frac{1}{i}$.
	Consider the following expression:
	\[
	   \lambda \binom{k}{t} v^t \sum_{i=0}^{\lambda-1} \binom{N}{i} p^i (1-p)^{N-i} \left(H_{N-i} - H_{N} + \log\left(\frac{1}{1-p}\right)\right).
	\]

	Since $H_k$ is an increasing function in $k$, $H_{N-i} - H_N \le 0$ for any $i \ge 0$.
	Then an estimate on the best number of rows $R$ in the initial stage can be found by determining when the following equation is satisfied:
	\[
	    \log\left(\frac{1}{1-p}\right)\lambda \binom{k}{t} v^t \sum_{i=0}^{\lambda-1} \binom{R}{i} p^i (1-p)^{R-i} < 1.
	\]
	Since $\log(\frac{1}{1-p})$ does not depend on $R$, we can apply the same methods and upper bounds as those used to prove \Cref{thm:standard_prob_analysis}.
 Following the same steps as that proof yields the theorem statement.
\end{proof}

\begin{figure}
    \centering
    \begin{tikzpicture}
\begin{axis}[no marks,
	legend pos=south east,
        legend style={at={(\legendxpos,\legendypos)},anchor=south east},
	width=0.9\textwidth,
        height=3in,
	xmin=\xmin,
	xmax=\xmax,
	xlabel={\xaxislabel},
	ylabel={\yaxislabel}]
\addplot[black,forget plot] table [mark=none, x=k, y=slj_no_sum_no_W, col sep=comma] {\filename};
\addplot[black,forget plot] table [mark=none, x=k, y=slj_no_sum_with_W, col sep=comma] {\filename};
\addplot[black,forget plot] table [mark=none, x=k, y=slj_with_sum, col sep=comma] {\filename};

\addplot[black,forget plot] table [mark=none, x=k, y=lll_no_sum_no_W, col sep=comma] {\filename};
\addplot[black,forget plot] table [mark=none, x=k, y=lll_no_sum_with_W, col sep=comma] {\filename};
\addplot[black,forget plot] table [mark=none, x=k, y=lll_with_sum, col sep=comma] {\filename};

\addplot[red] table [mark=none, x=k, y=two_stage_no_sum_no_W, col sep=comma] {\filename};
\addplot[blue] table [mark=none, x=k, y=two_stage_no_sum_with_W, col sep=comma] {\filename};
\addplot[orange] table [mark=none, x=k, y=two_stage_with_sum, col sep=comma] {\filename};

\legend{\Cref{thm:differentiation_general}, \Cref{thm:differentiation_general} + $W_{-1}$, \Cref{eqn:diff_general_setup_original}}
\end{axis}
\end{tikzpicture}
        \caption{Upper bounds for $\CAN_\lambda(t,k,v)$ when $t = 6, k \le 2000, v = 4$, and $\lambda = 12$ from \Cref{thm:differentiation_general}, the same type of bound but with $W_{-1}$, and \Cref{eqn:diff_general_setup_original}.
    These three bounds are the lowest lines within each of the three sections.
    The other upper bounds are from the three SLJ results from \Cref{fig:slj_bounds_v4_t6_l12} and the three LLL results from \Cref{fig:lll_bounds_v4_t6_l12}.}
    \label{fig:two_stage_bounds_v4_t6_l12}
\end{figure}
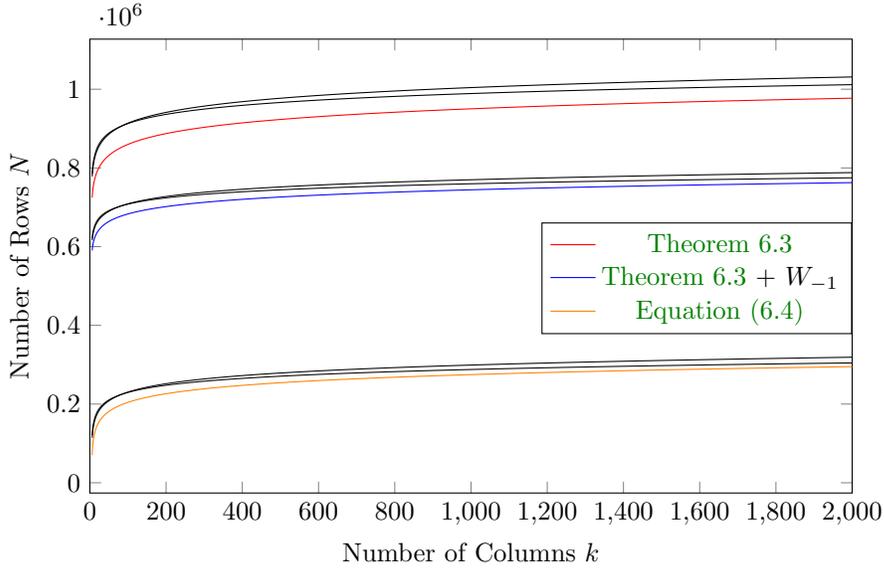

Jensen's inequality shows that $\log k + \gamma + \frac{1}{2k+1}\le H_k \le \log k + \gamma + \frac{1}{2k-1}$, where $\gamma$ is the Euler-Mascheroni constant.
Therefore, it is in principle possible to improve the constants in \Cref{thm:differentiation_general} further.
One such method would involve improving the estimate $H_{N-i} - H_N$ for all $0 \le i < \lambda$ by utilizing the upper and lower bounds produced by Jensen's inequality.
However, this improved estimate does not appear to give a closed-form solution in terms of $N$ and $W$ in the context of the entire equation.

For \Cref{thm:differentiation_2,thm:differentiation_general}, the number of rows added was $\lambda$ for each uncovered interaction, even if the interactions considered are covered some number of times already.
Naturally, determining the expected number of interactions covered $i$ times, for each $0 \le i < \lambda$, would be desired, since then we can improve the upper bound to $\lambda-i$ more rows for each of these interactions.
However, the event of whether an interaction is $i$-covered, and the similar event for the same interaction being $j$-covered, for $i \ne j$, are not independent. 
Thus one cannot simply determine the expectation for each and incorporate them in the failure probability. 

\subsection{Graph Coloring}

The methods of \Cref{thm:differentiation_2,thm:differentiation_general} construct a new row for every interaction that is not covered.
For each of these rows, $k-t$ entries can be arbitrarily set; if $k$ is much larger than $t$, then this strategy does not properly make use of such other entries. 
We outline a method that provides better bounds in general, which is a generalization of the methods in other work by Sarkar and Colbourn \cite{sarkar_two-stage_2019}.
First construct a random array as was done in the first stages of \Cref{thm:differentiation_2,thm:differentiation_general}; suppose that the interactions left uncovered are $I_1, \cdots, I_\rho$.
Further, suppose that for each such interaction $I_j$, it has been covered $d(I_j)$ times; note that $0 \le d(I_j) < \lambda$.
Construct a graph $G = (V, E)$ where each vertex in $V$ is a pair $(I, r)$, where $I$ is one of the uncovered interactions, and $d(I) + 1 \le r \le \lambda$.
We construct an edge between vertices $(I_1, r_1)$ and $(I_2, r_2)$ if $I_1 = I_2$ (and $r_1 \ne r_2$), or $I_1$ cannot be placed in the same row as $I_2$; we outline this criterion next. 
Call $I_1, I_2$ \emph{compatible} if for all columns $I_1$ and $I_2$ have in common, their entries are equal; $I_1, I_2$ are \emph{incompatible} otherwise. 
Then $\chi(G)$, the minimum number of colors to properly color $G$, is the minimum number of rows to cover all the remaining interactions. 
This generalizes Sarkar and Colbourn's result \cite{sarkar_two-stage_2019} in that it is extended to any $\lambda \ge 1$.

After $N$ rows are selected, \Cref{thm:differentiation_2,thm:differentiation_general} analyze the expected number of uncovered interactions.
For graph coloring, this is the same as the number of vertices in $G$; instead we analyze the expected number of edges in $G$.
Using an appropriate generalization of Sarkar and Colbourn's methods, this expectation $r$ is equal to:
\[
    r = \frac{1}{2} \binom{k}{t} v^t \sum_{i=1}^t \binom{t}{i} \binom{k-t}{t-i} (v^t - v^{t-i}) f(v^t) f(v^t - v^{t-i}),
\]
where
\[
    f(x) = \sum_{j=0}^{\lambda-1} \binom{N}{j} \frac{1}{x^j} \left(1 - \frac{1}{x}\right)^{N-j}.
\]
We use the standard upper bound of $\chi(G) \le \frac{1}{2} + \sqrt{2r + \frac{1}{4}}$ \cite{deistel_graph_2005} to give an upper bound on $\CAN_\lambda(t,k,v)$. 
One could hope to use the ideas from \Cref{thm:standard_prob_analysis} to give a strong analytical upper bound for a two-stage approach with graph coloring, but we were not able to do so. 
Instead, we performed a crude approximation and upper bounds using the methods in the proof of \Cref{thm:standard_prob_analysis} of the optimal choice  for the initial number of rows $N$ that yields a small overall number of rows, using the methods from \Cref{thm:differentiation_general}.
We found that one should choose approximately an $N$ that minimizes
\begin{equation}\label{eqn:two_stage_graph_analytical_approx}
    N + \sqrt{q}N^\lambda \left(1-\frac{1}{v^t}\right)^{N+1}
\end{equation}    
for $N \ge 0$, where
\[
q = \frac{\binom{k}{t} v^{2t} }{2} \frac{e^{2\lambda}}{\lambda^{2\lambda}} \frac{(k-2t)\binom{k-t}{t-1}(k! (k-2t)! - ((k-t)!)^2)}{t((k-t)!)^2}.
\]

To generate such an array algorithmically, we follow a similar procedure as in \Cref{sec:slj}: first calculate a value of $N$ that minimizes \Cref{eqn:two_stage_graph_analytical_approx} and create an array $A$ with $N$ rows (and $k$ columns) with entries chosen uniformly-at-random, independently.
If $A$ has at most $r$ interactions that are not $\lambda$-covered, then proceed to the second stage. 
Otherwise, go back to the first step.
In the second stage, generate the graph $G$ as described above using the interactions not $\lambda$-covered.
Generate a random vertex coloring using the upper bound on the chromatic number $\chi(G)$ of colors; if this is proper vertex coloring, then proceed.
Otherwise, try another coloring. 
For each of the color classes, form a row with the interactions within the color class put into that row (by definition, these interactions must be compatible).
Add all such rows to the partial array generated in the first stage.

\begin{figure}
    \centering
    \begin{tikzpicture}
\begin{axis}[no marks,
	legend pos=south east,
        legend style={at={(\legendxpos,\legendypos)},anchor=south east},
	width=0.9\textwidth,
        height=3in,
	xmin=\xmin,
	xmax=\xmax,
	xlabel={\xaxislabel},
	ylabel={\yaxislabel}]
\addplot[black,forget plot] table [mark=none, x=k, y=slj_no_sum_no_W, col sep=comma] {\filename};
\addplot[black,forget plot] table [mark=none, x=k, y=slj_no_sum_with_W, col sep=comma] {\filename};
\addplot[black,forget plot] table [mark=none, x=k, y=slj_with_sum, col sep=comma] {\filename};

\addplot[black,forget plot] table [mark=none, x=k, y=lll_no_sum_no_W, col sep=comma] {\filename};
\addplot[black,forget plot] table [mark=none, x=k, y=lll_no_sum_with_W, col sep=comma] {\filename};
\addplot[black,forget plot] table [mark=none, x=k, y=lll_with_sum, col sep=comma] {\filename};

\addplot[black,forget plot] table [mark=none, x=k, y=two_stage_no_sum_no_W, col sep=comma] {\filename};
\addplot[black,forget plot] table [mark=none, x=k, y=two_stage_no_sum_with_W, col sep=comma] {\filename};
\addplot[black,forget plot] table [mark=none, x=k, y=two_stage_with_sum, col sep=comma] {\filename};

\addplot[red] table [mark=none, x=k, y=two_stage_coloring, col sep=comma] {\filename};

\legend{Graph Coloring}
\end{axis}
\end{tikzpicture}
\caption{Upper bounds for $\CAN_\lambda(t,k,v)$ when $t = 6, k \le 2000, v = 4$, and $\lambda = 12$ from Graph Coloring in the second stage, which is the lowest line presented.
     The other upper bounds are from the three SLJ results from \Cref{fig:slj_bounds_v4_t6_l12}, the three LLL results from \Cref{fig:lll_bounds_v4_t6_l12}, and the three results from \Cref{fig:two_stage_bounds_v4_t6_l12}.}
\label{fig:graph_color_bounds_v4_t6_l12}
\end{figure}
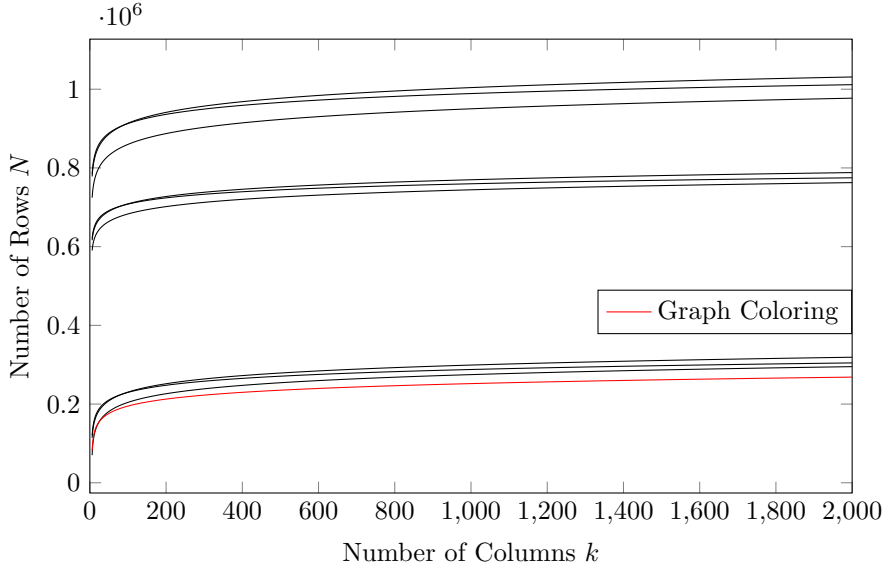

Unfortunately there does not appear to be an analytical solution directly to \Cref{eqn:two_stage_graph_analytical_approx}, including the use of $W$.
Nevertheless, as before we report computationally-found upper bounds from graph coloring (in the second stage) against the three SLJ bounds from \Cref{fig:slj_bounds_v4_t6_l12}, the three LLL bounds from \Cref{fig:lll_bounds_v4_t6_l12}, and the three bounds from \Cref{fig:two_stage_bounds_v4_t6_l12} for $\CAN_\lambda(t,k,v)$ where $t=6, k \le 2000, v=4$, and $\lambda=12$.
As is evident, graph coloring produces the smallest upper bounds as the number of rows in the second stage will be on average smaller than the corresponding number from the other two-stage results. 
The initial array produced in the first stage is random, so there is no information that can be inferred about the distribution of the number of edges, or other properties of the graph; this was noted by Sarkar and Colbourn \cite{sarkar_two-stage_2019}.


\section{Fixed Number of Rows}\label{sec:relaxation}

In this section we address the question of how large $\lambda$ can be when the other parameters $N, t, k, v$ are fixed.
We start with a straightforward application of the methods from \Cref{thm:standard_prob_analysis}.

\begin{theorem}\label{thm:lower_bound_lambda}
    Let $N, t, k, v, p$ be as in \Cref{thm:standard_prob_analysis}.
    Let 
    \[
        b = \binom{k}{t} v^t (1-p)^{N+1} \sqrt{\frac{1}{1-2p}}.
    \]
    Suppose that $N$ is sufficiently large, in that $1/e^N < b < 1$.
    Then there exists a $\CA_\lambda(N; t, k, v)$, where:
    \[
        \lambda = \left\lceil N \exp\left( 1 + W_{-1}\left(\frac{\log b}{eN}\right) \right) \right\rceil - 1.
    \]
\end{theorem}

\begin{proof}
    The proof starts the same way as \Cref{thm:standard_prob_analysis} through the calculation resulting in \Cref{eqn:standard_upper_bd}, reproduced here:
    \[
        \binom{k}{t} v^t (1-p)^{N-\lambda+1} a \left(\frac{eN}{\lambda}\right)^{\lambda} = 1.
    \]
    Recall that in the statement of \Cref{thm:standard_prob_analysis}, $a = \sqrt{\frac{(1-p)^{2\lambda} - p^{2\lambda}}{1-2p}}$.
    We upper bound $a$ by $\sqrt{\frac{(1-p)^{2\lambda}}{1-2p}}$, which cancels the $(1-p)^{-\lambda}$ term from \Cref{eqn:standard_upper_bd}.
    Therefore, if
    \[
        \binom{k}{t} v^t (1-p)^{N+1} \frac{1}{\sqrt{1-2p}} \left(\frac{eN}{\lambda}\right)^{\lambda} = 1,
    \]
    then any decrease in $\lambda$ will cause the left side of this equation to be strictly smaller than 1, which justifies the existence of a covering array.
    Substitution of $b$ into the above expression and solving for $\lambda$ yields:
    \[
        \lambda = N \exp\left(1 + W\left(\frac{\log b}{eN}\right)\right).
    \]
    If $1/e^N < b < 1$, then $-1/e < \frac{\log b}{eN} < 1$, and thus there are two real solutions. 
    We choose the negative branch $W_{-1}$ for the same reasons as in the proof of \Cref{thm:standard_prob_analysis}.
    Simple algebra yields the theorem statement.
\end{proof}

\begin{corollary}\label{cor:lower_bound_lambda_without_lambert}
    Let $N, k, v, t, b$ be as in \Cref{thm:lower_bound_lambda}.
    Then if 
    \[
        \lambda \le \left\lceil N \exp\left(-\frac{1 + e\log(-N/\log(b))}{e-1}\right) \right\rceil - 1,
    \]
    then there exists a $\CA_\lambda(N; t, k, v)$. 
\end{corollary}

\begin{proof}
    As was done in \Cref{cor:phrase_first_bound_without_lambert}, we use the lower bound on $W_{-1}$ by Alzharani and Salem \cite{alzahrani_sharp_2018} to the result of \Cref{thm:lower_bound_lambda}.
    Simple algebra yields the corollary statement.
\end{proof}

It is possible to improve upon these bounds using the Lov\'asz local lemma in a similar fashion to the proofs in \Cref{sec:LLL}.
Set $b' = e(\binom{k}{t} - \binom{k-t}{t}) v^t (1-p)^{N+1} \sqrt{1/(1-2p)}$.
Then there exists a $\CA_\lambda(N;t,k,v)$ where $\lambda = \lceil N \exp(1 + W_{-1}(\log b' / (eN))) \rceil$ by the methods in \Cref{thm:lower_bound_lambda}, and a similar quantity without $W_{-1}$ by the methods in \Cref{cor:lower_bound_lambda_without_lambert}.
As was done in the proof of \Cref{thm:standard_lll}, since the local lemma does not require the left-hand side being strictly less than 1, we are not required to subtract 1 in the bounds for $\lambda$ here. 

\section{Conclusion}\label{sec:conclusion}

In this paper we generalized previous work on covering array upper bounds for when the redundancy $\lambda$ is at least 1.
We utilized the Lambert $W$ function, specifically the negative branch $W_{-1}$ and lower bounds for that, to derive asymptotically optimal upper bounds for $\CAN_\lambda(t,k,v)$.
This framework was applied to a Stein-Lov\'asz-Johnson approach using a simple probabilistic model, then to the Lov\'asz local lemma, and then to a two-stage approach.
It is clear from our results that the two-stage approach produces the best upper bounds in general, but only can offer improvements in the constants as the SLJ results were asymptotically optimal when $v, t$ are fixed.

We outline some ideas for future work; obvious avenues involve improving the upper bounds even further.
Some possibilities are with group actions \cite{sarkar_upper_2017}, covering perfect hash families \cite{sarkar_covering_2016}, and alternate probability distributions \cite{godbole_t-covering_1996}. 
Some of the bounds they produce are the best-known in general for $\CA_1$'s; does the same hold true in the higher-index setting? 
Can a graph decomposition technique improve the second stage even further?

Further, what lower bounds can be determined? 
\Cref{cor:ca_tight_bound_first} shows that the best lower bound must be of the form $c_1 \log k + c_2 \lambda$ for some constants $c_1, c_2$ that are independent of $k, \lambda$ but may depend on $v, t$. 
There remains a large gap between the upper and lower bounds for $c_2$.
Lower bounds for covering arrays are not well explored, but the methods by Choi et al. \cite{choi_structures_2012} may be generalizable to $\lambda > 1$.

As noted by Sarkar and Colbourn \cite{sarkar_two-stage_2019}, merits of applying more than two stages would potentially be fruitful but need to have justification.
Dougherty \cite{dougherty_genetic_2020} gives empirical evidence that more stages produce better bounds through the use of a genetic algorithm; however, analyzing analytically what algorithm for each stage is best, and how much the algorithm should contribute to eventually having every interaction covered at least $\lambda$ times, would likely be very difficult.
One possible direction is to use an exact method, such as constraint satisfaction \cite{hnich_constraint_2006}, to find the optimal number of rows in the second stage. 
A generalization of the methods in this paper would have to be applied to so-called ``variable-strength'' covering arrays \cite{raaphorst_variable_2018}; this generalization has logarithmic growth in the number of interactions for fixed $\lambda, v, t$.
Therefore, the number of rows in the second stage would be logarithmic in the number of interactions not-yet $\lambda$-covered, whereas the coloring approach in \Cref{sec:two_stage} produces approximately a square root number of additional rows. 
Algorithmic generation of the second stage could use the methods by Dougherty et al. \cite{dougherty_algorithmic_2022}.
However, the trade-off in rows would be an increase in computation time (and storage, depending on the method).

Better, known methods in the first stage can inform the second stage how to proceed.
For example, the Kleitman and Spencer \cite{kleitman_families_1973} and Katona \cite{katona_two_1973} approach for determining $\CAN_1(2,k,2)$ involve finding the maximum number of subsets of $\{1, \cdots, N\}$ with pairwise nonempty intersection and symmetric difference. 
The optimal construction involves selecting all subsets of size approximately half the number of rows $N$.
Analyzing this method shows that only a very small set of pairs of interactions are covered exactly once. 
Additionally, are there recursive constructions for covering arrays that make use of higher $\lambda$ (instead of simply being a generalization of existing recursive constructions)?

Another potential avenue is for variable-strength covering arrays \cite{moura_lovasz_2018} or partial covering arrays \cite{sarkar_partial_2018}, in which not all sets of $t$ columns need to cover all $v^t$ interactions within them, and other suitable generalizations.
If $N$ and $\lambda$ are fixed, how many $t$ sets of columns can have their interactions covered?

\section{Acknowledgments}

The views expressed in this article are those of the author(s) and do not reflect the official policy or position of the Department of the Army, Department of Defense, or the U.S. Government.



\bibliographystyle{siamplain}
\bibliography{ryan_references}
\end{document}

%% file: ex_shared.tex
\usepackage{lipsum}
\usepackage{amsfonts}
\usepackage{graphicx}
\usepackage{epstopdf}
\usepackage{algorithmic}
\ifpdf
  \DeclareGraphicsExtensions{.eps,.pdf,.png,.jpg}
\else
  \DeclareGraphicsExtensions{.eps}
\fi

\newsiamremark{remark}{Remark}
\newsiamremark{hypothesis}{Hypothesis}
\crefname{hypothesis}{Hypothesis}{Hypotheses}
\newsiamthm{claim}{Claim}

\usepackage{tikz}
\usepackage{pgfplots}
\pgfplotsset{compat=1.13}

\headers{Upper Bounds for Covering Arrays of Higher Index}{M. R. Calbert and R. E. Dougherty}

\title{Upper Bounds for Covering Arrays of Higher Index}

\author{Mason R. Calbert\thanks{Department of Electrical Engineering \& Computer Science, United States Military Academy, West Point, NY
  (\email{mason.calbert@westpoint.edu}).}
\and Ryan E. Dougherty\thanks{Department of Electrical Engineering \& Computer Science, United States Military Academy, West Point, NY 
  (\email{ryan.dougherty@westpoint.edu}).}
}

\usepackage{amsopn}

\makeatletter
\newcommand*{\addFileDependency}[1]{
  \typeout{(#1)}
  \@addtofilelist{#1}
  \IfFileExists{#1}{}{\typeout{No file #1.}}
}
\makeatother

\newcommand*{\myexternaldocument}[1]{%
    \externaldocument{#1}%
    \addFileDependency{#1.tex}%
    \addFileDependency{#1.aux}%
}